\def\draft{n}

\documentclass[12pt]{amsart}
\usepackage{fullpage,amssymb,epic,eepic,amscd,epsfig}

\theoremstyle{plain}

\newtheorem{theorem}{Theorem}
\newtheorem{proposition}{Proposition}[section]
\newtheorem{lemma}[proposition]{Lemma}
\newtheorem{corollary}[proposition]{Corollary}

\theoremstyle{definition}

\newtheorem{question}[proposition]{Question}

\theoremstyle{remark}

\def\printname#1{
	\if\draft y
		\smash{\makebox[0pt]{\hspace{-0.5in}
			\raisebox{8pt}{\tt\tiny #1}}}
	\fi
}

\newcommand{\psdraw}[2]
        {\begin{array}{c} \hspace{-1.3mm}
	\raisebox{-4pt}{\psfig{figure=draws/#1.ps,width=#2}}
	\hspace{-1.9mm}\end{array}}
\newlength{\standardunitlength}
\catcode`\@=11
\long\def\@makecaption#1#2{%
    \vskip 10pt
    
\setbox\@tempboxa\hbox{
      \small\sf{\bfcaptionfont #1. }\ignorespaces #2}%
    \ifdim \wd\@tempboxa >\captionwidth {%
        \rightskip=\@captionmargin\leftskip=\@captionmargin
        \unhbox\@tempboxa\par}%
      \else
        \hbox to\hsize{\hfil\box\@tempboxa\hfil}%
    \fi}
\font\bfcaptionfont=cmssbx10 scaled \magstephalf
\newdimen\@captionmargin\@captionmargin=2\parindent
\newdimen\captionwidth\captionwidth=\hsize
\catcode`\@=12

\def\lbl#1{\label{#1}\printname{#1}}

\def\BZ{\mathbb Z}
\def\BQ{\mathbb Q}

\def\A{\mathcal A}
\def\B{\mathcal B}

\def\M{\mathcal M}

\def\La{\Lambda}

\def\ihs{integral homology 3-sphere}  
\def\fti{finite type invariant}       

\def\la{\langle}
\def\ra{\rangle}

\def\s{\sigma}
\def\we{\wedge}
\def\bwe{\bigwedge}

\def\a{\alpha}

\def\Ga{\Gamma}

\def\BBI{\mathbb I}
\def\BBJ{\mathbb J}

\def\la{\langle}
\def\ra{\rangle}
\def\we{\wedge}
\def\Ga{\Gamma}
\def\La{\Lambda}
\def\sp{\mathfrak s\mathfrak p}
\def\tilCphi{\tilde{\mathcal C}(\phi)}
\def\Cphi{\mathcal C(\phi)}

\def\Vert{\text{\,Vert}}
\def\Edge{\text{\,Edge}}
\def\Flag{\text{\,Flag}}
\def\Sym{\text{Sym}}
\def\Im{\text{Im}}
\def\sgn{\text{ sgn}}
\def\ot{\otimes}
\def\Sp{\text{Sp}}
\def\Gr{\text{Gr}}
\def\BBU{\mathbb U}

\begin{document}

\title[Some $IHX$ relations  and symplectic
representation theory]{Some $IHX$-type relations on trivalent graphs
and symplectic representation theory}

\author{Stavros Garoufalidis}
\address{Department of Mathematics \\
         Brandeis University \\
         Waltham, MA 02254-9110, U.S.A. }
\email{stavros@oscar.math.brandeis.edu}
\thanks{The  first author was partially supported by NSF grant 
       DMS-95-05105. \newline
       This and related preprints can also be obtained at
{\tt http://www.math.brown.edu/$\sim$stavrosg }}

\author{Hiroaki Nakamura}
\address{Department of Mathematics\\
        Tokyo Metropolitan University \\
        Tokyo 192-0397, Japan}
\email{h-naka@comp.metro-u.ac.jp}

\date{
This edition: May 12, 1998
\hspace{0.5cm} First edition: April 5, 1997 
}

\begin{abstract}
We consider two types of graded algebras 
(with graded actions by the symplectic Lie algebra)
that arise in the study of the mapping class group, 
and describe their symplectic invariants in terms of 
algebras on trivalent graphs. 
\end{abstract}

\maketitle


\section{Introduction}
\lbl{sec.intro}

Let $\sp_g$ be the Lie algebra of symplectic matrices of degree $2g$
over the rational numbers. 
In recent studies related to the structure of 
the {\em surface mapping class group},
several authors \cite{KM,Mo,Ha,HL} have encountered
a certain distinguished quotient $\B$ of the exterior algebra $\La U$, 
where $U$ is an irreducible $\sp_g$-module
isomorphic to $\Lambda^3H/H$, $H$ is the fundamental $\sp_g$-module
and $\La^k$ is the $k^{th}$ exterior functor.
The second exterior component $\La^2 U$ is 
decomposed as an $\sp_g$-module in the following way:
\begin{equation}
\lbl{eq.symp}
\Lambda^2U\cong [1^6]_{\sp} \oplus [1^4]_{\sp} \oplus
[1^2]_{\sp} \oplus [0]_{\sp} \oplus [2^21^2]_{\sp} \oplus [2^2]_{\sp}
\quad (g\ge 6), 
\end{equation}
where $[\lambda]_{\sp}$ denotes an irreducible $\sp_g$-module 
corresponding to a partition $\lambda$, and 
the algebra $\B$ mentioned above is defined by
\begin{equation}
\B:=\La U /([2^2]_{\sp}),
\end{equation}
where  $[2^2]_{\sp}\subset\La^2U$ according to the decomposition 
\eqref{eq.symp}.  

In fact, the algebra $\B$ appears
in the Hodge theoretic study of the mapping class group $M_g$ 
(of a closed genus $g$ surface) by R. Hain \cite{Ha}, 
who established a theory of mixed Hodge structure for
the Torelli group $T_g:=\ker(M_g\to\Sp_g(\BZ))$.
Introducing the unipotent kernel  $\mathfrak u_g$ of the 
`relative Malcev completion' of the map $M_g\to \Sp_g(\mathbb Q)$,
Hain showed that the universal envelope $\BBU\Gr^W\mathfrak u_g$ 
of the weight graded Lie algebra $\Gr^W\mathfrak u_g$ is
quadratic dual to the lowest weight subalgebra of the 
continuous cohomology $H_{cts}^\ast(\mathfrak u_g)$.
Then, using the fact $\Gr^W\mathfrak u_g$ is generated by 
the weight $-1$ component isomorphic to $U$, 
he gave a natural isomorphism
$$
\B \cong \bigoplus_{k\ge 1}W_kH_{cts}^k(\mathfrak u_g).
$$
See \cite[Prop. 9.9]{HL}.
>From this, the $\sp_g$-invariants of $\B$ are
naturally mapped into the cohomology algebra $H^{\ast}(M_g,\BQ)$.

On the other hand, S.Morita \cite{Mo} introduced a Weyl-type \cite{W} 
interpretation of $(\La U)^{\sp}$, the space of $\sp_g$-invariants 
of $\La U$, by the algebra of trivalent graphs $\Cphi$ 
(see 2.1 below). Then, using the 
generalized 
Morita-Mumford-Miller classes (with twisted coefficients)
from \cite{K},
N.Kawazumi and S.Morita \cite{KM,Mo} related explicitly
the trivalent graphs `lying in' $(\La U)^{\sp}$ with the 
Morita-Mumford-Miller classes $e_n\in H^{2n}(M_g,\BQ)$,
and showed that $(\La U)^{\sp}$ surjects onto the subalgebra 
$\BQ[e_1,e_2,\dots]$ in $H^{\ast}(M_g,\BQ)$.

Motivated by the above mentioned works \cite{Ha,HL,K,KM,Mo}
(as well as by the stable independence of the classes $e_i$'s 
due to Miller, Morita),
one would like to understand stably
the structure of the $\sp_g$-invariant algebra
$\B^{\sp}=(\La U /([2^2]_{\sp}))^{\sp}$
purely combinatorially in 
terms of trivalent graphs and their relations.

Let $\Cphi$ denote the commutative
graded algebra generated by the trivalent graphs that (possibly)
contain  multiple
edges and 1-loops (where a 1-loop is an edge which begins and 
ends on the same vertex). By definition, the 
multiplication in $\Cphi$ is given by disjoint union of graphs,
and the degree of a trivalent graph is half the number of vertices.
Let also $loop$ denote the ideal generated by
graphs containing a 1-loop. 

Define $IH_0$ to be the ideal of $\Cphi$
`identifying' $I=H$ 
(with 4 distinct edges connected to a central edge)
locally in trivalent graphs. 

\begin{theorem}
\lbl{thm.new1}
There exists a stable isomorphism of graded algebras
$$
\Cphi/(IH_0, loop) 
\buildrel{\sim}\over\longrightarrow
(\La U /([2^2]_{\sp}))^{\sp}
$$
which multiplies degrees by $2$. 
Here, ``stable'' means that for each degree $m$, the 
map is an isomorphism for $g \geq 3m$.
\end{theorem}

It is easy to see that $\Cphi/(IH_0, loop)$ is a polynomial algebra 
(freely) generated by the $IH_0$-classes of connected trivalent graphs
without 1-loops,
i.e., having one generator in every positive degree (cf. 3.3 (c) below).
See also \cite[p. 639]{KM} for a relevant discussion involving 
Kontsevich's primitive factors ``$H_0(\text{Out}F_n,\BQ)$'' 
\cite{Ko3}.

In recent studies related to finite type 3-manifold invariants 
\cite{GO,LMO,L}, perturbative Chern-Simons theory \cite{Wi, Ka, Ko4, RW}
and graph cohomology \cite{Ko1, Ko2}, 
the algebra  $\mathcal A(\phi) = 
\tilCphi/(AS,IHX)$ plays a crucial role, where 
$\tilCphi$ is the algebra generated by vertex-oriented trivalent
graphs and  $AS,IHX$ are the relations shown in
figure \ref{ASIHX}.

\begin{figure}[htpb]
$$ \printname{ASIHX}
	\setlength{\unitlength}{0.03\standardunitlength}
	\begin{array}{c}  \hspace{-1.7mm}
        	\raisebox{-8pt}{\begingroup\makeatletter\ifx\SetFigFont\undefined
\def\x#1#2#3#4#5#6#7\relax{\def\x{#1#2#3#4#5#6}}%
\expandafter\x\fmtname xxxxxx\relax \def\y{splain}%
\ifx\x\y   
\gdef\SetFigFont#1#2#3{%
  \ifnum #1<17\tiny\else \ifnum #1<20\small\else
  \ifnum #1<24\normalsize\else \ifnum #1<29\large\else
  \ifnum #1<34\Large\else \ifnum #1<41\LARGE\else
     \huge\fi\fi\fi\fi\fi\fi
  \csname #3\endcsname}%
\else
\gdef\SetFigFont#1#2#3{\begingroup
  \count@#1\relax \ifnum 25<\count@\count@25\fi
  \def\x{\endgroup\@setsize\SetFigFont{#2pt}}%
  \expandafter\x
    \csname \romannumeral\the\count@ pt\expandafter\endcsname
    \csname @\romannumeral\the\count@ pt\endcsname
  \csname #3\endcsname}%
\fi
\fi\endgroup
\begin{picture}(9396,1157)(0,-10)
\thicklines
\path(6238.154,632.231)(6312.000,533.000)(6293.538,655.308)
\put(6462.000,595.500){\arc{325.000}{2.7468}{6.6780}}
\path(8430.462,655.308)(8412.000,533.000)(8485.846,632.231)
\put(8262.000,595.500){\arc{325.000}{2.7468}{6.6780}}
\path(238.154,1007.231)(312.000,908.000)(293.538,1030.308)
\put(462.000,970.500){\arc{325.000}{2.7468}{6.6780}}
\path(1588.154,557.231)(1662.000,458.000)(1643.538,580.308)
\put(1812.000,520.500){\arc{325.000}{2.7468}{6.6780}}
\path(3613.154,332.231)(3687.000,233.000)(3668.538,355.308)
\put(3837.000,295.500){\arc{325.000}{2.7468}{6.6780}}
\path(238.154,107.231)(312.000,8.000)(293.538,130.308)
\put(462.000,70.500){\arc{325.000}{2.7468}{6.6780}}
\path(2488.154,557.231)(2562.000,458.000)(2543.538,580.308)
\put(2712.000,520.500){\arc{325.000}{2.7468}{6.6780}}
\path(4063.154,332.231)(4137.000,233.000)(4118.538,355.308)
\put(4287.000,295.500){\arc{325.000}{2.7468}{6.6780}}
\path(12,983)(912,983)
\path(12,83)(912,83)
\path(1812,983)(1812,83)
\path(2712,983)(2712,83)
\path(3612,983)(4512,83)
\path(3612,83)(3987,458)
\path(4137,608)(4512,983)
\path(6462,533)(6462,83)
\path(7812,983)(8262,533)(8712,983)
\path(8262,533)(8262,83)
\path(6012,983)(6462,533)(6912,983)
\path(462,983)(462,83)
\path(1812,533)(2712,533)
\path(3837,308)(4287,308)
\put(1137,458){\makebox(0,0)[lb]{$=$}}
\put(3162,458){\makebox(0,0)[lb]{$-$}}
\put(7212,458){\makebox(0,0)[lb]{$+$}}
\put(9012,458){\makebox(0,0)[lb]{$=0$}}
\end{picture} }
        	\hspace{-1.9mm}
	\end{array}
 $$
\caption{The $AS$ and $IHX$ relation on vertex-oriented trivalent 
graphs.}\lbl{ASIHX}
\end{figure}

In fact, letting $\M$ be the vector space spanned by the isomorphism
classes of (oriented) integral homology 3-spheres,
and $\widehat{\A}(\phi)$ the completion of $\A(\phi)$
(with respect to the graph degrees), 
Le-Murakami-Ohtsuki \cite{LMO}
constructed a map $\Omega:\M\to \widehat{\A}(\phi)$ which turns out 
\cite{L,GO} to be 
the `universal' finite type invariant of integral homology 3-spheres,
in the sense that $\Omega$ induces an isomorphism 
between  the graded space of $\M$ (with respect to the Ohtsuki filtration)
and  ${\A}(\phi)$.
The proof of Theorem \ref{thm.new1} provides an alternative 
characterization of $\A(\phi)$ as follows. 
Noting the $\sp$-decomposition of $\La^2(\Sym^3 H)$
(where $\Sym^n$ denotes the $n$-th symmetric tensor functor)   
\begin{equation}
\La^2(\Sym^3 H) = 
[5,1]_{\sp}\oplus [4]_{\sp}\oplus [3^2]_{\sp}\oplus
[2^2]_{\sp}\oplus [1^2]_{\sp}\oplus [0]_{\sp}
\quad (g\ge 2),
\end{equation}
\begin{theorem} 
\lbl{thm.new2}
Stably, we have an isomorphism of graded algebras 
(which multiplies degrees by $2$)
$$ 
\mathcal A(\phi)
\buildrel{\sim}\over\longrightarrow
 (\La(\Sym^3 H)/([4]_{\sp}))^{\sp}.
$$
\end{theorem}
In the above theorem  ``stably'' means that for each degree $m$, 
the map is an isomorphism for $g \geq m$.

Since $[4]_{\sp}\cong\Sym^4 H$, this result
is closely related to the study of
the Lie algebra $Ham=\bigoplus_{m\ge 2}\Sym^m H$ 
of formal Hamiltonian vector fields with Poisson brackets
\cite{Ko2,Ko3}. In fact, the kernel  $Ham^1$ of 
$Ham\to \Sym^2H\cong\sp_g$ is generated by 
$\Sym^3 H$, thus inducing a natural homomorphism
$\La(\Sym^3 H)/([4]_{\sp})\to H^\ast(Ham^1,\BQ)$
(hence also $\mathcal A(\phi)\to H^\ast(Ham,\BQ)$).
In view of comparison of two situations of 
our above theorems, it seems an interesting (future) subject to determine 
the kernel/cokernels of the latter homomorphisms.

The idea of the present paper arose during conversations 
of the first named author with S. Morita in July 1996 
and during a visit of the second named author to Harvard 
in the fall of 1996.
We wish to thank T. Gocho and  P. Vogel for useful conversations 
and especially N.Kawazumi and S. Morita for enlightening 
communications and crucial remarks.

\section{A reduction of Theorem \ref{thm.new1}}
\lbl{sec.2}

\subsection{Trivalent graphs and $\sp$-invariant tensors}
\lbl{sub.21}

In the course of proving Theorems \ref{thm.new1} and \ref{thm.new2},
we will deduce several other results of independent interest, which
we first formulate.

As was mentioned in the introduction, S. Morita 
\cite{Mo} (using the device of H. Weyl [W]) introduced  maps
\begin{equation}
\lbl{eq.va}
\Cphi \to (\La \La^3 H)^{\sp} \quad\text{ and }\quad
\Cphi/(loop) \to (\La U )^{\sp},
\end{equation}
and showed that they are
(stable) isomorphisms of graded algebras (which multiply
degrees by $2$). ``Stable'' here means that in each degree $m$, the above
maps are isomorphisms for $g \geq 3m$.
These maps essentially place  a copy of the symplectic form 
on each edge of a trivalent graph, but for the calculations needed 
below, we prefer to give coordinatewise definitions of them.

Let $\Ga$ be a {\em trivalent graph} with {\em vertex} set $\Vert(\Ga)$
and {\em edge} set $\Edge(\Ga)$ and let $\Flag(\Ga)$ be 
the set of {\em flags}, where a flag is by definition a  pair
 consisting of a vertex and an incident half-edge.  
Since each vertex has three adjacent flags, the cardinality
$|\Flag(\Ga)|$ is equal to $3|\Vert(\Ga)|=2|\Edge(\Ga)|$.
We call $m=1/6|\Flag(\Ga)|$ the degree of the trivalent 
graph $\Ga$. 
A {\em total ordering} $\tau$ of $\Ga$ consists, by definition,
of the following data:
\begin{itemize}
\item
a linear ordering of vertices: 
$\Vert(\Ga)=\{v_1,\dots,v_{2m}\}$;
\item
 a linear ordering of $\Flag(v)=\{f_1(v),f_2(v),f_3(v)\}$
for each $v\in \Vert(\Ga)$;
\item
an ordering of $\Flag(e)=\{f_+(e),f_-(e)\}$ for each 
$e\in\Edge(\Ga)$.
\end{itemize}
Such a $\tau$ is called  $\we$-admissible if
it satisfies the condition:
$$
\sgn
\pmatrix 
f_1(v_1) & f_2(v_1) & f_3(v_1) &f_1(v_2) & ...f_3(v_{2m}) \\
f_+(e_1) & f_-(e_1) & f_+(e_2) &f_-(e_2) & ...f_-(e_{3m})
\endpmatrix
=1
$$
for every linear ordering of edges: 
$\Edge(\Ga)=\{e_1,\dots,e_{3m}\}$.
Note that the sign is not changed under edge permutations.

Let $H$ be an $\sp_g$-module equipped with a standard symplectic basis
$\{x_1,\dots,x_g,y_1,\dots,y_g\}$
with $\la x_i,y_j\ra=\delta_{ij}=-\la y_i,x_j\ra$,
$\la x_i,x_j\ra=\la y_i,y_j\ra=0$ over $\mathbb Q$. 
Given a totally ordered trivalent graph $(\Ga,\tau)$ of degree $m$, 
we define 
$f_{3i+j}:=f_j(v_i)$ ($0\le i<2m$, $j=1,2,3$)
and $OR:=\{f_+(e)\}_{e\in\Edge(\Ga)}$, and shall
define an $\sp$-invariant 
$\alpha_{(\Ga,\tau)}\in H^{\otimes 6m}$ 
as follows.
First, regard $H^{\otimes 6m}$ as the linear
combinations of the words in $\{x_i,y_i\}_{i=1}^g$ 
of length $6m$, and write $x_{-i}=y_i$, $x_0=0$.
Let $I$ be the set of indices
$\mathbf i=(i_1,\dots,i_{6m}) \in \{-g,-g+1,\ldots,g-1,g \}^{6m}$
such that $i_k+i_l=0$ if and only if $f_k$ and $f_l$ share
an edge. 
Then, we define
\begin{equation*}
\alpha_{(\Ga,\tau)}=\sum_{\mathbf i\in I} \sgn(\mathbf i)x_{\mathbf i}
\quad\in (H^{\otimes 6m})^{\sp},
\end{equation*}
where $\sgn(\mathbf i)=\prod_{f_k\in OR}\sgn(i_k)$ and
$x_{\mathbf i}=x_{i_1}\cdots x_{i_{6m}}$.
The image of $\alpha_{(\Ga,\tau)}$ via the standard projection 
$H^{\otimes 6m}\to \Lambda^{2m} \Lambda^3 H $ is {\em independent} 
of $\tau$ (as long as $\tau$ is chosen to be $\we$-admissible), 
and will be denoted by 
$\alpha_\Ga\in(\Lambda^{2m} \Lambda^3H )^{\sp}$.
Since the kernel of $\Lambda^3H\to U$ equals to 
$H\we \sum_i(x_i\we y_i)$, it is easy to see that 
the image of $\a_\Ga$ in $\La^{2m} U $
(denoted by the same symbol) vanishes exactly when
the graph $\Ga$ has a 1-loop.
Extending the map $\Ga\mapsto \alpha_\Ga$ linearly, 
we obtain the (stable) isomorphisms 
in \eqref{eq.va}. 

Let us now  introduce our basic $\sp$-homomorphisms:
$f_I$, $f_H$, $f_X:H^{\otimes 4}\to \Lambda^2 \Lambda^3 H $
by setting the images 
of $t=t_1 \ot t_2 \ot t_3 \ot t_4 \in H^{\ot 4}$
as follows: 
\begin{eqnarray*}
f_I(t)  & = &  \sum_{i=1}^g
(t_1 \we t_2 \we x_i) \we  (t_3 \we t_4 \we y_i) -
(t_1 \we t_2 \we y_i) \we  (t_3 \we t_4 \we x_i), \\ 
f_{H}(t)  & = &  \sum_{i=1}^g 
(t_1 \we t_3 \we x_i) \we  (t_4 \we t_2 \we y_i) - 
(t_1 \we t_3 \we y_i) \we  (t_4 \we t_2 \we x_i), \\ 
f_X(t)  & = &  \sum_{i=1}^g 
(t_1 \we t_4 \we x_i) \we  (t_2 \we t_3 \we y_i) - 
(t_1 \we t_4 \we y_i) \we   (t_2 \we t_3 \we x_i).
\end{eqnarray*}
(We also use the same symbols to denote the compositions
of $f_I,f_H,f_X$ with the projection 
$\Lambda^2 \Lambda^3H\to\Lambda^2 U$ respectively.)

We now define, for scalars $a,b,c$, the $\sp_g$-homomorphism:
$$ f_{a,b,c}= a f_I + b f_H + c f_X,$$
and we  denote the  composite of $f_{a,b,c}$
with the projection $\La^2 \La^3 H\to\La^2 U$ by $\bar{f}_{a,b,c}$.
Furthermore, given an embedding of graphs $I \hookrightarrow \Ga$, 
one may associate three graphs $\Ga=\Ga_I,\Ga_H,\Ga_X$
by replacing adjacent relations of 4 edges into $I$
by those into $I$, $H$, $X$ as indicated in figure \ref{gIHX}
respectively. 
Using this notation, for any triple of scalars $(a,b,c)$, 
we define $I_{a,b,c}\subset\Cphi$ to be
the ideal generated by the $a\Ga_I+b\Ga_H+c\Ga_X$
for all pairs $I \hookrightarrow \Ga$.

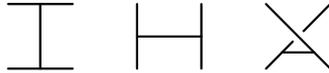
\begin{figure}[htpb]
$$ \printname{gIHX}
	\setlength{\unitlength}{0.03\standardunitlength}
	\begin{array}{c}  \hspace{-1.7mm}
        	\raisebox{-8pt}{\begingroup\makeatletter\ifx\SetFigFont\undefined
\def\x#1#2#3#4#5#6#7\relax{\def\x{#1#2#3#4#5#6}}%
\expandafter\x\fmtname xxxxxx\relax \def\y{splain}%
\ifx\x\y   
\gdef\SetFigFont#1#2#3{%
  \ifnum #1<17\tiny\else \ifnum #1<20\small\else
  \ifnum #1<24\normalsize\else \ifnum #1<29\large\else
  \ifnum #1<34\Large\else \ifnum #1<41\LARGE\else
     \huge\fi\fi\fi\fi\fi\fi
  \csname #3\endcsname}%
\else
\gdef\SetFigFont#1#2#3{\begingroup
  \count@#1\relax \ifnum 25<\count@\count@25\fi
  \def\x{\endgroup\@setsize\SetFigFont{#2pt}}%
  \expandafter\x
    \csname \romannumeral\the\count@ pt\expandafter\endcsname
    \csname @\romannumeral\the\count@ pt\endcsname
  \csname #3\endcsname}%
\fi
\fi\endgroup
\begin{picture}(4524,939)(0,-10)
\thicklines
\path(12,912)(912,912)
\path(12,12)(912,12)
\path(1812,912)(1812,12)
\path(2712,912)(2712,12)
\path(3612,912)(4512,12)
\path(3612,12)(3987,387)
\path(4137,537)(4512,912)
\path(1812,462)(2712,462)
\path(3837,237)(4287,237)
\path(462,912)(462,12)
\end{picture} }
        	\hspace{-1.9mm}
	\end{array}
 $$
\caption{The graphs $I, H, X$.}\lbl{gIHX}
\end{figure}

In the next section, we will show the following
\begin{proposition}
\lbl{thm.A}
The stable isomorphism of equation \eqref{eq.va} induces 
stable isomorphisms of graded algebras 
\begin{eqnarray*}
\Cphi/I_{a,b,c} \cong
(\Lambda \Lambda^3 H/(\Im\, f_{a,b,c}))^{\sp} \\
\Cphi/(I_{a,b,c}\cup loop) \cong
(\Lambda U/(\Im\,\bar{f}_{a,b,c}))^{\sp}
\end{eqnarray*}
which multiply degrees by $2$.
\end{proposition}
\noindent
Letting $IHX$ (resp. $IH$) denote the ideal $I_{1,1,1}$ 
(resp. $I_{1,-1,0}$),
and  $IH_0$ denote the ideal generated by
 the restricted form of the $IH$-relation 
where $I,H$ are connected by 4 distinct edges, we obtain the following:
\begin{corollary}
\lbl{cor.rel}
For every pair  $(\BBI,\BBJ)$  of  the following table: 
\newline
\vspace{-0.3cm}
\begin{center}
\begin{tabular}{|r|c|} \hline
$ \BBI $ & $ \BBJ \subset\La^2U$ \\  \hline
$IHX$ & $[1^4]_{\sp} \oplus
[1^2]_{\sp} \oplus [0]_{\sp} $ \\ \hline
$IH$  & $[2^2]_{\sp}\oplus [1^2]_{\sp}\oplus [0]_{\sp} $ \\ \hline
$IH_0 $ & $ [2^2]_{\sp} $ \\ \hline
\end{tabular}
\end{center}
we have a stable isomorphism: 
$$ 
\Cphi/(\BBI \cup loop)  \cong (\Lambda U/ (\BBJ))^{\sp}.
$$
\end{corollary}

We will also show the following
\begin{proposition}
\lbl{prop.rel}
(a) $(\La U /([1^4]_{\sp} \oplus[1^2]_{\sp} \oplus [0]_{\sp}))^{\sp}$
vanishes in positive degrees up to degree $12$. \newline
(b) $(\La U /([2^2]_{\sp}\oplus [1^2]_{\sp}\oplus [0]_{\sp}))^{\sp}$
vanishes in positive degrees. \newline
(c)
$(\La U /([2^2]_{\sp}))^{\sp}$ is a free polynomial algebra having 
one generator in every even positive degree.
\end{proposition}

Theorem \ref{thm.new1} follows immediately from Corollary \ref{cor.rel}
and Proposition \ref{prop.rel} above.

\section{Proofs for the results of Section 2}
\subsection{Proof of Proposition \ref{thm.A}}

Let $J$ be the ideal of $\La \La^3 H$ generated by the image
of $f_{a,b,c}$ and $J_{2m}$ its homogeneous part of degree $2m$.
It suffices to show that $J_{2m}^{\sp}$ is generated by the 
$a\alpha_{\Ga_I}+b\alpha_{\Ga_H}+ c \alpha_{\Ga_X}$
for all pairs $I \hookrightarrow \Ga$ with $deg(\Ga)=m$.
Since $\La \La^3 H $ is skew-commutative, $J_{2m}=
( \Im\, f_{a,b,c})  \we ( \La^{2m-2 }\La^3 H) $.
We thus have a surjection
$F_{a,b,c}:H^{\ot 4}\ot H^{\ot 6m-6}\to J_{2m}$ 
factoring through $f_{a,b,c}\ot id^{\ot 6m-6}$.
If we define 
$F_I,F_H,F_X:H^{\ot 4}\ot H^{\ot 6m-6}\to \La^{2m} \La^3 H$ 
by replacing $f_{a,b,c}$ by $f_I,f_H,f_X$ in the above, then
obviously we have $F_{a,b,c}=aF_I+bF_H+cF_X$.
Let $C:H^{\ot 2}\to \mathbb Q$ denote the canonical contraction
which maps $\sum x_i\ot y_i-y_i\ot x_i$ to 1, and define
$C\ot F_\ast:H^{\ot 6m}\to \La^{2m} \La^3 H$ 
$(\ast=\{a,b,c\},I,H,X)$
in such a way that the domain components of $C$ (resp. 
$f_\ast$ part of $F_\ast$)
share the third and sixth (resp. first, second, 
fourth and fifth) positions of $H^{\ot 6m}$.
Still $C\ot F_{a,b,c}$ gives a surjection onto $J_{2m}$, and
the semisimplicity of $\sp$-representations implies that
$J_{2m}^{\sp}$ is generated by the images of $\sp$-invariants 
$\alpha_{(\Ga,\tau)}$ of $H^{\ot 6m}$ via $C\ot F_{a,b,c}$,
where $(\Ga,\tau)$ runs over trivalent graphs of degree $m$
with $\we$-admissible
total orderings such that $f_3=f_+(e')$, $f_6=f_-(e')$ 
for some $e'\in\Edge(\Ga)$.
Given such a $(\Ga,\tau)$, 
construct three trivalent graphs
$\Ga=\Ga_I, \Ga_H, \Ga_X$ by replacing 
$I$ by $H$ and $X$ for the latter two, and
give their total orderings $\tau_I,\tau_H,\tau_X$
so that $\tau_I=\tau$ and $\tau_H,\tau_X$ differ
from $\tau$ only locally in the parts `$H,X$'
as indicated in figure \ref{IHXwe}.

\begin{figure}[htpb]
$$ \printname{IHXwe}
	\setlength{\unitlength}{0.03\standardunitlength}
	\begin{array}{c}  \hspace{-1.7mm}
        	\raisebox{-8pt}{\begingroup\makeatletter\ifx\SetFigFont\undefined
\def\x#1#2#3#4#5#6#7\relax{\def\x{#1#2#3#4#5#6}}%
\expandafter\x\fmtname xxxxxx\relax \def\y{splain}%
\ifx\x\y   
\gdef\SetFigFont#1#2#3{%
  \ifnum #1<17\tiny\else \ifnum #1<20\small\else
  \ifnum #1<24\normalsize\else \ifnum #1<29\large\else
  \ifnum #1<34\Large\else \ifnum #1<41\LARGE\else
     \huge\fi\fi\fi\fi\fi\fi
  \csname #3\endcsname}%
\else
\gdef\SetFigFont#1#2#3{\begingroup
  \count@#1\relax \ifnum 25<\count@\count@25\fi
  \def\x{\endgroup\@setsize\SetFigFont{#2pt}}%
  \expandafter\x
    \csname \romannumeral\the\count@ pt\expandafter\endcsname
    \csname @\romannumeral\the\count@ pt\endcsname
  \csname #3\endcsname}%
\fi
\fi\endgroup
\begin{picture}(7286,1755)(0,-10)
\thicklines
\path(4218.462,905.308)(4200.000,783.000)(4273.846,882.231)
\put(4050.000,845.500){\arc{325.000}{2.7468}{6.6780}}
\path(5851.154,657.231)(5925.000,558.000)(5906.538,680.308)
\put(6075.000,620.500){\arc{325.000}{2.7468}{6.6780}}
\path(6693.462,680.308)(6675.000,558.000)(6748.846,657.231)
\put(6525.000,620.500){\arc{325.000}{2.7468}{6.6780}}
\path(676.154,1332.231)(750.000,1233.000)(731.538,1355.308)
\put(900.000,1295.500){\arc{325.000}{2.7468}{6.6780}}
\path(1068.462,455.308)(1050.000,333.000)(1123.846,432.231)
\put(900.000,395.500){\arc{325.000}{2.7468}{6.6780}}
\path(3318.462,905.308)(3300.000,783.000)(3373.846,882.231)
\put(3150.000,845.500){\arc{325.000}{2.7468}{6.6780}}
\path(3150,1308)(3150,408)
\path(4050,1308)(4050,408)
\path(3270.000,888.000)(3150.000,858.000)(3270.000,828.000)
\path(3150,858)(4050,858)
\path(5850,1308)(6750,408)
\path(5850,408)(6225,783)
\path(6375,933)(6750,1308)
\path(6075,633)(6525,633)
\path(6405.000,603.000)(6525.000,633.000)(6405.000,663.000)
\path(450,1308)(1350,1308)
\path(450,408)(1350,408)
\path(900,1308)(900,408)
\path(870.000,528.000)(900.000,408.000)(930.000,528.000)
\put(0,1383){\makebox(0,0)[lb]{$e_2$}}
\put(1425,1383){\makebox(0,0)[lb]{$e_1$}}
\put(0,183){\makebox(0,0)[lb]{$e_4$}}
\put(4425,783){\makebox(0,0)[lb]{$1$}}
\put(2625,783){\makebox(0,0)[lb]{$2$}}
\put(6000,108){\makebox(0,0)[lb]{$1$}}
\put(6450,108){\makebox(0,0)[lb]{$2$}}
\put(825,1608){\makebox(0,0)[lb]{$1$}}
\put(825,33){\makebox(0,0)[lb]{$2$}}
\put(4125,1383){\makebox(0,0)[lb]{$e_1$}}
\put(6825,1383){\makebox(0,0)[lb]{$e_1$}}
\put(2700,1383){\makebox(0,0)[lb]{$e_2$}}
\put(5400,1383){\makebox(0,0)[lb]{$e_2$}}
\put(1425,183){\makebox(0,0)[lb]{$e_3$}}
\put(2700,183){\makebox(0,0)[lb]{$e_4$}}
\put(5400,183){\makebox(0,0)[lb]{$e_4$}}
\put(4125,183){\makebox(0,0)[lb]{$e_3$}}
\put(6825,183){\makebox(0,0)[lb]{$e_3$}}
\end{picture} }
        	\hspace{-1.9mm}
	\end{array}
 $$
\caption{Three graphs $I,H,X$ together with
their total ordering, where shown are particular orderings 
of the trivalent vertices together with flag orientations 
for the respective vertices and internal edges.
The external edges $e_1,...,e_4$ are assumed
(flag-)oriented in the same way on the three graphs.}\lbl{IHXwe}
\end{figure}
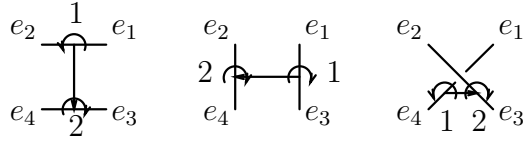

Then, it is easy to see that 
$C\ot F_\ast(\alpha_{(\Ga,\tau)})=\alpha_{\Ga_\ast}$
for $\ast=I,H,X$ and hence that 
$C\ot F_{a,b,c}(\alpha_{(\Ga,\tau)})
=a\alpha_{\Ga_I}+b\alpha_{\Ga_H}+c\alpha_{\Ga_X}$.
Conversely, given three trivalent graphs $\Ga_I,\Ga_H,\Ga_X$ 
which coincide except in their distinguished parts $I,H,X$
respectively, we may associate $\we$-admissible total 
orderings $\tau_I,\tau_H,\tau_X$ 
on the three graphs $\Ga_I,\Ga_H,\Ga_X$
such that they coincide with each other
outside the `$IHX$-parts' and such that their orderings
inside the `$IHX$-parts' are as indicated in figure \ref{IHXwe}
with the middle edges shared by $f_3, f_6$ on the three
graphs.
It is not difficult then to see that 
$C\ot F_{a,b,c}(\alpha_{(\Ga_I,\tau_I)})
=a\alpha_{\Ga_I}+b\alpha_{\Ga_H}+c\alpha_{\Ga_X}$.
This completes the proof of Proposition \ref{thm.A}.
\qed

\subsection{Proof of Corollary \ref{cor.rel}}

(i) The case $\BBI=IHX$: It is easy to see that $f_{IHX}$
(i.e., $f_{1,1,1}$) is alternating, 
hence factors through 
$\Lambda^4 H\cong[1^4]_{\sp} \oplus[1^2]_{\sp} \oplus [0]_{\sp}.$ 
We can choose highest weight vectors for the $\sp$-decomposition of 
$\La^4 H$ as follows: $x_1 \we x_2 \we x_3 \we x_4$ for $[1^4]_{\sp}$,
$\sum_{i=1}^g x_i \we y_i \we x_1 \we x_2 $ for $[1^2]_{\sp}$
and $\sum_{i,j=1}^n x_i \we y_i \we x_j \we y_j $ for $[0]_{\sp}$.
A computation of their images in $\La^2 U$ under $f_{IHX}$
shows that $\Lambda^4 H$ is embedded into $\La^2 U$.
For example, with the temporary abbreviation of $\alpha_\Gamma$ by $\Gamma$
(due to typesetting reasons), we have
\begin{eqnarray*}
f_{IHX}(\sum_{i=1}^n  x_i \ot y_i \ot x_1 \ot x_2 ) & 
= &  \printname{circle2}
	\setlength{\unitlength}{0.03\standardunitlength}
	\begin{array}{c}  \hspace{-1.7mm}
        	\raisebox{-8pt}{\begingroup\makeatletter\ifx\SetFigFont\undefined
\def\x#1#2#3#4#5#6#7\relax{\def\x{#1#2#3#4#5#6}}%
\expandafter\x\fmtname xxxxxx\relax \def\y{splain}%
\ifx\x\y   
\gdef\SetFigFont#1#2#3{%
  \ifnum #1<17\tiny\else \ifnum #1<20\small\else
  \ifnum #1<24\normalsize\else \ifnum #1<29\large\else
  \ifnum #1<34\Large\else \ifnum #1<41\LARGE\else
     \huge\fi\fi\fi\fi\fi\fi
  \csname #3\endcsname}%
\else
\gdef\SetFigFont#1#2#3{\begingroup
  \count@#1\relax \ifnum 25<\count@\count@25\fi
  \def\x{\endgroup\@setsize\SetFigFont{#2pt}}%
  \expandafter\x
    \csname \romannumeral\the\count@ pt\expandafter\endcsname
    \csname @\romannumeral\the\count@ pt\endcsname
  \csname #3\endcsname}%
\fi
\fi\endgroup
\begin{picture}(5561,1632)(0,-10)
\thicklines
\put(825.000,1158.000){\arc{900.000}{1.5708}{4.7124}}
\path(945.000,1638.000)(825.000,1608.000)(945.000,1578.000)
\put(825.000,1158.000){\arc{900.000}{4.7124}{7.8540}}
\put(2972.727,1103.455){\arc{1013.595}{2.9521}{4.6181}}
\put(4772.728,1103.455){\arc{1013.595}{2.9521}{4.6181}}
\path(3047.292,1626.566)(2925.000,1608.000)(3041.641,1566.832)
\put(2877.273,1103.455){\arc{1013.595}{4.8067}{6.4727}}
\path(4847.292,1626.566)(4725.000,1608.000)(4841.641,1566.832)
\put(4677.272,1103.455){\arc{1013.597}{4.8067}{6.4727}}
\path(375,408)(1275,408)
\path(825,708)(825,408)
\path(795.000,528.000)(825.000,408.000)(855.000,528.000)
\path(2475,1008)(2475,408)
\path(3375,1008)(3375,408)
\path(4275,408)(5175,1008)
\path(4275,1008)(4650,783)
\path(4800,633)(5175,408)
\path(2595.000,1038.000)(2475.000,1008.000)(2595.000,978.000)
\path(2475,1008)(3375,1008)
\path(4425,483)(5025,483)
\path(4905.000,453.000)(5025.000,483.000)(4905.000,513.000)
\put(750,933){\makebox(0,0)[lb]{$1$}}
\put(0,108){\makebox(0,0)[lb]{$x_2$}}
\put(1350,108){\makebox(0,0)[lb]{$x_1$}}
\put(675,33){\makebox(0,0)[lb]{$2$}}
\put(1725,933){\makebox(0,0)[lb]{$+$}}
\put(3750,933){\makebox(0,0)[lb]{$+$}}
\put(2250,108){\makebox(0,0)[lb]{$x_2$}}
\put(5100,108){\makebox(0,0)[lb]{$x_2$}}
\put(3225,108){\makebox(0,0)[lb]{$x_1$}}
\put(4125,108){\makebox(0,0)[lb]{$x_1$}}
\put(2625,708){\makebox(0,0)[lb]{$1$}}
\put(3075,708){\makebox(0,0)[lb]{$2$}}
\put(4200,558){\makebox(0,0)[lb]{$1$}}
\put(5100,558){\makebox(0,0)[lb]{$2$}}
\end{picture} }
        	\hspace{-1.9mm}
	\end{array}
 \\
& = & \printname{circle3}
	\setlength{\unitlength}{0.03\standardunitlength}
	\begin{array}{c}  \hspace{-1.7mm}
        	\raisebox{-8pt}{
\begingroup\makeatletter\ifx\SetFigFont\undefined
\def\x#1#2#3#4#5#6#7\relax{\def\x{#1#2#3#4#5#6}}%
\expandafter\x\fmtname xxxxxx\relax \def\y{splain}%
\ifx\x\y   
\gdef\SetFigFont#1#2#3{%
  \ifnum #1<17\tiny\else \ifnum #1<20\small\else
  \ifnum #1<24\normalsize\else \ifnum #1<29\large\else
  \ifnum #1<34\Large\else \ifnum #1<41\LARGE\else
     \huge\fi\fi\fi\fi\fi\fi
  \csname #3\endcsname}%
\else
\gdef\SetFigFont#1#2#3{\begingroup
  \count@#1\relax \ifnum 25<\count@\count@25\fi
  \def\x{\endgroup\@setsize\SetFigFont{#2pt}}%
  \expandafter\x
    \csname \romannumeral\the\count@ pt\expandafter\endcsname
    \csname @\romannumeral\the\count@ pt\endcsname
  \csname #3\endcsname}%
\fi
\fi\endgroup
\begin{picture}(3686,1632)(0,-10)
\thicklines
\put(825.000,1158.000){\arc{900.000}{1.5708}{4.7124}}
\path(945.000,1638.000)(825.000,1608.000)(945.000,1578.000)
\put(825.000,1158.000){\arc{900.000}{4.7124}{7.8540}}
\put(2972.727,1103.455){\arc{1013.595}{2.9521}{4.6181}}
\path(3047.292,1626.566)(2925.000,1608.000)(3041.641,1566.832)
\put(2877.273,1103.455){\arc{1013.595}{4.8067}{6.4727}}
\path(375,408)(1275,408)
\path(825,708)(825,408)
\path(795.000,528.000)(825.000,408.000)(855.000,528.000)
\path(2475,1008)(2475,408)
\path(3375,1008)(3375,408)
\path(2595.000,1038.000)(2475.000,1008.000)(2595.000,978.000)
\path(2475,1008)(3375,1008)
\put(750,933){\makebox(0,0)[lb]{$1$}}
\put(0,108){\makebox(0,0)[lb]{$x_2$}}
\put(1350,108){\makebox(0,0)[lb]{$x_1$}}
\put(675,33){\makebox(0,0)[lb]{$2$}}
\put(2250,108){\makebox(0,0)[lb]{$x_2$}}
\put(3225,108){\makebox(0,0)[lb]{$x_1$}}
\put(1725,1008){\makebox(0,0)[lb]{$+2$}}
\put(2625,708){\makebox(0,0)[lb]{$1$}}
\put(3000,708){\makebox(0,0)[lb]{$2$}}
\end{picture}
 }
        	\hspace{-1.9mm}
	\end{array}

\end{eqnarray*}
in $\Lambda^2 \Lambda^3H $ whose first (resp. second)
term vanishes (resp. remains) when projected to $\Lambda^2 U$.

(ii) The case of $\BBI=IH$: 
Let $f_{IH}=f_I-f_H=f_{1,-1,0}$.
Since $f_{IH}(t_1\ot t_2\ot t_3\ot t_4)$ 
is invariant under
the variable changes of $t_1\leftrightarrow t_4$, 
$t_2\leftrightarrow t_3$ and $t_1,t_4\leftrightarrow t_2,t_3$
respectively, we find that $f_{IH}$ factors through 
$\Sym^2(\Sym^2H)\cong [4]_{\sp}\oplus [2^2]_{\sp}\oplus 
[1^2]_{\sp}\oplus [0]_{\sp}$.
Since the target space does not have $[4]_{\sp}$, $f_{IH}$
is zero on this component.
The highest weight vectors of the other components of 
$\Sym^2(\Sym^2H)$ can be taken as 
$(x_1x_2)(x_1x_2)-(x_1^2)(x_2^2)$ for $[2^2]_{\sp}$,
$\sum_{j=1}^g(x_1x_j)(x_2y_j)-(x_1y_j)(x_2x_j)$ for $[1^2]_{\sp}$
and $\sum_{k,l=1}^g(x_kx_l)(y_ky_l)-(x_ky_l)(x_ly_k)$ for
$[0]_{\sp}$.
Mapping these vectors by $f_{IH}$, we find that 
$[2^2]_{\sp}$, $[1^2]_{\sp}$ and $[0]_{\sp}$ remain
nontrivally in $\Lambda^2 U$.

(iii) The case $\BBI=IH_0$:
We have to re-examine the proof of Proposition \ref{thm.A} carefully.
Observe first that the $IH$-relation can be classified into
the three types $IH_0$, $IH_1$, $IH_2$ according to 
the number of connected edges to `$I$-graph' being 4,3 or 2.

\begin{figure}[htpb]
$$ \printname{eqIH}
	\setlength{\unitlength}{0.03\standardunitlength}
	\begin{array}{c}  \hspace{-1.7mm}
        	\raisebox{-8pt}{\begingroup\makeatletter\ifx\SetFigFont\undefined
\def\x#1#2#3#4#5#6#7\relax{\def\x{#1#2#3#4#5#6}}%
\expandafter\x\fmtname xxxxxx\relax \def\y{splain}%
\ifx\x\y   
\gdef\SetFigFont#1#2#3{%
  \ifnum #1<17\tiny\else \ifnum #1<20\small\else
  \ifnum #1<24\normalsize\else \ifnum #1<29\large\else
  \ifnum #1<34\Large\else \ifnum #1<41\LARGE\else
     \huge\fi\fi\fi\fi\fi\fi
  \csname #3\endcsname}%
\else
\gdef\SetFigFont#1#2#3{\begingroup
  \count@#1\relax \ifnum 25<\count@\count@25\fi
  \def\x{\endgroup\@setsize\SetFigFont{#2pt}}%
  \expandafter\x
    \csname \romannumeral\the\count@ pt\expandafter\endcsname
    \csname @\romannumeral\the\count@ pt\endcsname
  \csname #3\endcsname}%
\fi
\fi\endgroup
\begin{picture}(14646,944)(0,-10)
\thicklines
\path(12,914)(912,914)
\path(12,14)(912,14)
\path(462,914)(462,14)
\path(1812,914)(1812,14)
\path(2712,914)(2712,14)
\path(1812,464)(2712,464)
\put(1137,389){\makebox(0,0)[lb]{$=$}}
\put(4742.357,464.000){\arc{910.714}{1.4173}{4.8659}}
\path(4812,914)(5712,914)
\path(4812,14)(5712,14)
\path(5262,914)(5262,14)
\path(7212,914)(7212,14)
\path(8112,914)(8112,14)
\path(7212,464)(8112,464)
\put(7142.357,464.000){\arc{910.714}{1.4173}{4.8659}}
\put(6012,389){\makebox(0,0)[lb]{$=$}}
\path(13212,914)(13212,14)
\path(14112,914)(14112,14)
\path(13212,464)(14112,464)
\put(10142.357,464.000){\arc{910.714}{1.4173}{4.8659}}
\put(11182.000,464.000){\arc{910.824}{4.5581}{8.0083}}
\put(13142.357,464.000){\arc{910.714}{1.4173}{4.8659}}
\put(14182.000,464.000){\arc{910.824}{4.5581}{8.0083}}
\path(10212,914)(11112,914)
\path(10212,14)(11112,14)
\path(10662,914)(10662,14)
\put(12012,389){\makebox(0,0)[lb]{$=$}}
\end{picture} }
        	\hspace{-1.9mm}
	\end{array}
 $$
\caption{The three possible forms of the $IH$ relation
in trivalent graphs (with no orderings).
On the left, the four external edges of $I$ are 
assumed distinct.}\lbl{eqIH}
\end{figure}
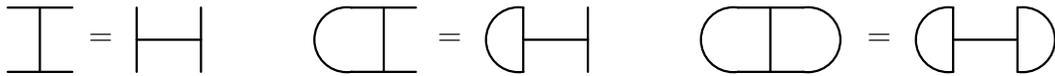 
\noindent
Replace $f_{1,-1,0}=f_I-f_H$ in the proof of Theorem 1 
by its composite $f'$ with the projection to the $[2^2]_{\sp}$-part,
and set $J=\bigoplus J_m$ to be the ideal generated by the
$\Im (f')$.
Form a surjection $C\ot F':H^{\ot 6m}\to J_{2m}$
factoring through $C\ot f'\ot id^{\ot 6m-6}$.
Then,  we have to show  that the collection of 
the images of $C\ot F'(\alpha_{(\Ga,\tau)})$ in $\La^{2m} U $
(where $(\Ga,\tau)$ runs over all trivalent graphs of degree $m$
with $\we$-admissible total orderings) 
coincides exactly with the collection of 
the $\alpha_{\Ga_I}-\alpha_{\Ga_H}$
where $(\Ga_I,\Ga_H)$ runs over all $IH_0$-related pairs
of trivalent graphs.
But the kernel of $\Im(f_{IH})\to[2^2]_{\sp}$ is
isomorphic to $\La^2 H=[1^2]_{\sp}\oplus[0]_{\sp}$
whose image in $\La^2 U $ is generated 
by the elements of the form
$$
\sum_{i,j}\sgn(ij)(x_k\we x_i\we x_j)
\we (x_l\we x_{-i}\we x_{-j}) 
\quad (k\ne l).
$$
{}From this it follows that, among the collection of the
$C\ot F_{IH}(\alpha_{(\Ga,\tau)})$, those that remain under
$C\ot F'$ are exactly the $\alpha_{\Ga_I}-\alpha_{\Ga_H}$
coming from $IH_0$-related pairs, as desired.

\subsection{Proof of Proposition \ref{prop.rel}}

(b): Using the $IH$ relation, we easily see that every (trivalent)
graph can be deformed so as to have a 1-loop, hence 
that $\Cphi/(IH \cup loop)$ vanishes
in positive degrees.

(c): Since $IH_0$ defines an equivalence relation on 
trivalent graphs, it follows that  $\Cphi/(IH_0 \cup loop)$ is a 
polynomial algebra on the graded set of 
equivalence classes of connected trivalent graphs without 
1-loops modulo the $IH_0$ relation.
By induction, it is easy to see that every connected graph of
degree $n$ without 1-loops is equivalent (modulo the $IH_0$ relation) to
the graph $E_n$ shown in Figure \ref{en}. 
Thus $\Cphi/(IH_0 \cup loop)$ is a polynomial algebra
on the classes of $E_n$ $(n\ge 1)$. 

\begin{figure}[htpb]
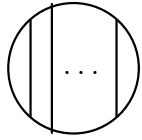

$$ \printname{en}
	\setlength{\unitlength}{0.03\standardunitlength}
	\begin{array}{c}  \hspace{-1.7mm}
        	\raisebox{-8pt}{\input draws/en.tex }
        	\hspace{-1.9mm}
	\end{array}
 $$
\caption{The graph $E_n$ of degree $n$, with $n$ vertical chords
on the circle.}\lbl{en}
\end{figure}

For (a), let us begin by collecting some elementary relations 
in $\Cphi/(IHX \cup loop)$.
We understand all graphs shown below as parts of trivalent graphs.
A cycle of length $n$ in a graph will be called an $n$-wheel, 
if it has $n$ distinct trivalent vertices.

\begin{figure}[htpb]
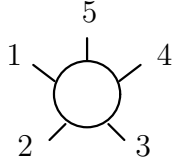

$$ \printname{wheel}
	\setlength{\unitlength}{0.03\standardunitlength}
	\begin{array}{c}  \hspace{-1.7mm}
        	\raisebox{-8pt}{\input draws/wheel.tex }
        	\hspace{-1.9mm}
	\end{array}
 $$
\caption{A wheel with five ordered legs.}\lbl{wheel}
\end{figure}

Given an $n$-wheel $w_I$ (where $I$ is the set of its $n$ external legs)
 and a permutation $\sigma$ of the symmetric group $\text{Sym}_I$,
we define another $n$-wheel $w_{\sigma(I)}$ to be the wheel whose
external legs are permuted by $\sigma$.

\begin{lemma}
\lbl{lem.n}
For a permutation $\s \in \text{Sym}_I$, and an $n$-wheel $w_I$,
we have:
\begin{equation}
w_{\s I} = \text{sgn}(\s) w_I \text{ modulo wheels with less than $n$
legs}
\end{equation}
\end{lemma}

\begin{proof}
Apply the $IHX$ relation, regarding $I$ as an 
arc of a wheel $w_I$ (i.e., an edge between
to adjacent trivalent vertices).
The wheel  appears in two
terms (with a different ordering of its legs) and an $n-1$ wheel
appears in the third term.
\end{proof}

\begin{corollary}
\lbl{cor.n}
A $2n$-wheel $w_{I}$ can be written as a linear combination 
of ($2n-1$)-wheels.
Indeed, observe that the sign of the permutation 
$\s=(1,2, \ldots, 2n)$
is $-1$, and that the graphs $w_{\s I}$ and  $w_{I}$ are isomorphic.
Furthermore, a $3$-wheel $w_J$ can be written as a linear combination
of $2$-wheels.
Indeed, observe that the graphs $w_J$ and $w_{(12)J}$ are isomorphic. 
Thus we conclude that $n$-wheels vanish in 
$\Cphi/(IHX \cup loop)$ for $n=1,2,3,4$.
\end{corollary}

\begin{lemma}
\lbl{lem.2p}
If a graph contains two pentagon cycles $C_5^1, C_5^2$ which intersect
at two consecutive edges, then it vanishes in $\Cphi/(IHX \cup loop)$.
\end{lemma}
\begin{proof}
Apply the $IHX$ relation, regarding $I$ as one of the two
consecutive common edges. We get that a sum of three terms vanish.
One of them is the graph in consideration, and the two others
contain squares, and thus by the above corollary vanish.
\end{proof}

\begin{corollary}
\lbl{cor.2p}
If a graph contains a pentagon cycle, then it
 can be written as
a linear combination of graphs that have at least 
$14$ vertices.
\end{corollary}

\begin{figure}[htpb]
$$ \psdraw{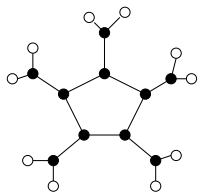}{1.0in} $$
\caption{A pentagon with $10$ vertices with $10$ external legs.}
\lbl{3tree2}
\end{figure}

\begin{proof}
By Corollary \ref{cor.n}, we have only to consider a graph without 
$(\le 4)$-wheels but with a pentagon.
Figure \ref{3tree2} shows that such a graph has at least
$10$ (black) vertices near the pentagon. 
Furthermore, there should appear $10$ (white) vertices 
adjacent to 5 of them. 
Using Lemma \ref{lem.2p}, we may assume that the white vertices 
are all distinct from the black ones,
but they need not be distinct from each other. Instead, at most $3$
of the white vertices can coincide with each other
(because this is part of a trivalent graph),  thus there are
at least $4$ distinct white vertices. Thus $\Ga$ has at least
$14$ vertices. 
\end{proof}

The {\em girth} $g(\Ga)$ of a (connected) graph $\Ga$ is, by definition, 
the minimum length of a cycle 
in $\Ga$ 
(\cite{B}; recall that a cycle is a closed path of distinct edges).
We now have the following Lemma 
(for a proof, see \cite[Theorem 1.2, p.105]{B}). 
\begin{lemma}
\lbl{lem.girth}
A trivalent graph with girth at least $g$, contains at least
$ 2^{(g+3)/2} -2$ (resp. $ 3 \cdot 2^{g/2} -2$) vertices if $g$ is odd 
(resp. $g$ is even).
\end{lemma}

Now, we  give the proof of (a) of Proposition \ref{prop.rel}.
Consider a connected trivalent graph $\Ga$.
If $\Ga$ has at most $8$ vertices, by Lemma \ref{lem.girth}, it contains
an $n$-gon for some $n=1,2,3,4$, and thus by corollary \ref{cor.n}
vanishes.
If $\Ga$ has girth exactly $5$, then by corollary \ref{cor.2p},
it has at least $14$ vertices.
If $\Ga$ has girth at least $6$, then by Lemma \ref{lem.girth},
it has at least $22$  vertices.
To sum up, if $\Ga$ has at most $12$ vertices, it follows that
$\Ga = 0 \in \Cphi/(IHX \cup loop)$. 
Thus, $\Cphi/(IHX \cup loop)$ vanishes 
in degrees $n=1,2,3,4,5,6$, and Corollary 
\ref{cor.rel} concludes our proof.

\begin{question}
Is the algebra $\Cphi/(IHX \cup loop)$ trivial in positive degrees?
\end{question}

\section{Proof of Theorem \ref{thm.new2}}

Let $\Ga$ be a {\it vertex oriented trivalent graph} i.e., 
a trivalent graph such that for each vertex $v\in \Edge(\Ga)$,
$\Flag(v)$ is given a cyclic ordering.
A {\it total ordering} $\tau$ of $\Ga$ is defined in the same 
way as in section \ref{sub.21} with an extra condition that 
the linear ordering of $\Flag(v)$ has the same sign as the 
given cyclic ordering on it.
We say such a $\tau$ to be $\Sym$-admissible if  
it satisfies the condition:
$$
\sgn
\pmatrix 
f_1(v_1) & f_2(v_1) & f_3(v_1) &f_1(v_2) & ...f_3(v_{2m}) \\
f_+(e_1) & f_-(e_1) & f_+(e_2) &f_-(e_2) & ...f_-(e_{3m})
\endpmatrix
=1
$$
for every linear ordering of edges: 
$\Edge(\Ga)=\{e_1,\dots,e_{3m}\}$.

\begin{figure}[htpb]
$$ \printname{IHXsym}
	\setlength{\unitlength}{0.03\standardunitlength}
	\begin{array}{c}  \hspace{-1.7mm}
        	\raisebox{-8pt}{\begingroup\makeatletter\ifx\SetFigFont\undefined
\def\x#1#2#3#4#5#6#7\relax{\def\x{#1#2#3#4#5#6}}%
\expandafter\x\fmtname xxxxxx\relax \def\y{splain}%
\ifx\x\y   
\gdef\SetFigFont#1#2#3{%
  \ifnum #1<17\tiny\else \ifnum #1<20\small\else
  \ifnum #1<24\normalsize\else \ifnum #1<29\large\else
  \ifnum #1<34\Large\else \ifnum #1<41\LARGE\else
     \huge\fi\fi\fi\fi\fi\fi
  \csname #3\endcsname}%
\else
\gdef\SetFigFont#1#2#3{\begingroup
  \count@#1\relax \ifnum 25<\count@\count@25\fi
  \def\x{\endgroup\@setsize\SetFigFont{#2pt}}%
  \expandafter\x
    \csname \romannumeral\the\count@ pt\expandafter\endcsname
    \csname @\romannumeral\the\count@ pt\endcsname
  \csname #3\endcsname}%
\fi
\fi\endgroup
\begin{picture}(7286,1410)(0,-10)
\thicklines
\path(5851.154,513.231)(5925.000,414.000)(5906.538,536.308)
\put(6075.000,476.500){\arc{325.000}{2.7468}{6.6780}}
\path(676.154,1188.231)(750.000,1089.000)(731.538,1211.308)
\put(900.000,1151.500){\arc{325.000}{2.7468}{6.6780}}
\path(676.154,288.231)(750.000,189.000)(731.538,311.308)
\put(900.000,251.500){\arc{325.000}{2.7468}{6.6780}}
\path(3826.154,738.231)(3900.000,639.000)(3881.538,761.308)
\put(4050.000,701.500){\arc{325.000}{2.7468}{6.6780}}
\path(6301.154,513.231)(6375.000,414.000)(6356.538,536.308)
\put(6525.000,476.500){\arc{325.000}{2.7468}{6.6780}}
\path(2926.154,738.231)(3000.000,639.000)(2981.538,761.308)
\put(3150.000,701.500){\arc{325.000}{2.7468}{6.6780}}
\path(3150,1164)(3150,264)
\path(4050,1164)(4050,264)
\path(3270.000,744.000)(3150.000,714.000)(3270.000,684.000)
\path(3150,714)(4050,714)
\path(5850,1164)(6750,264)
\path(5850,264)(6225,639)
\path(6375,789)(6750,1164)
\path(6075,489)(6525,489)
\path(6405.000,459.000)(6525.000,489.000)(6405.000,519.000)
\path(450,1164)(1350,1164)
\path(450,264)(1350,264)
\path(900,1164)(900,264)
\path(870.000,384.000)(900.000,264.000)(930.000,384.000)
\put(0,1239){\makebox(0,0)[lb]{$e_2$}}
\put(1425,1239){\makebox(0,0)[lb]{$e_1$}}
\put(4125,1239){\makebox(0,0)[lb]{$e_1$}}
\put(6825,1239){\makebox(0,0)[lb]{$e_1$}}
\put(2700,1239){\makebox(0,0)[lb]{$e_2$}}
\put(5400,1239){\makebox(0,0)[lb]{$e_2$}}
\put(0,39){\makebox(0,0)[lb]{$e_3$}}
\put(2700,39){\makebox(0,0)[lb]{$e_3$}}
\put(5400,39){\makebox(0,0)[lb]{$e_3$}}
\put(1425,39){\makebox(0,0)[lb]{$e_4$}}
\put(6825,39){\makebox(0,0)[lb]{$e_4$}}
\put(4125,39){\makebox(0,0)[lb]{$e_4$}}
\end{picture} }
        	\hspace{-1.9mm}
	\end{array}
 $$
\caption{Three graphs $I,H,X$ together with
their total ordering, where shown are particular orderings 
of the trivalent vertices together with flag orientations 
for the respective vertices and internal edges.
The external edges $e_1,...,e_4$ are assumed
(flag-)oriented in the same way on the three graphs.}\lbl{IHXsym}
\end{figure}
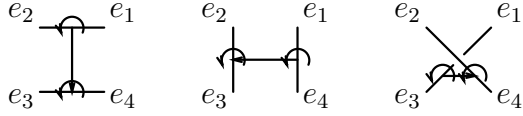

Given a totally ordered, vertex oriented trivalent graph
$(\Ga,\tau)$ of degree $m$, we can define an $\sp$-invariant
$\alpha_{(\Ga,\tau)}\in (H^{\otimes 6m})^{\sp}$
in the same way as in section \ref{sub.21}.
Moreover we see that its projection image in $\La^{2m}\Sym^3H$
is independent of $\tau$ as long as $\tau$ is chosen to be
$\Sym$-admissible. We denote this well-defined image by $\alpha_\Ga$.
The mapping $\alpha:\tilCphi\to (\La\Sym^3 H)^{\sp}$ 
($\Ga\mapsto \alpha_\Ga$)
turns out to be factoring through the AS-relation, inducing the
stable isomorphism $\tilCphi/(AS)\cong (\La\Sym^3 H)^{\sp}$.

Now we shall introduce our ``oriented version'' of the
basic $\sp$-homomorphisms:
$f_I$, $f_H$, $f_X:H^{\otimes 4}\to \Lambda^2 \Sym^3 H $
by setting the images 
of $t=t_1 \ot t_2 \ot t_3 \ot t_4 \in H^{\ot 4}$
as follows:   
\begin{eqnarray*}
f_I(t)  & = &  
\sum_{i=1}^g
(t_1t_2x_i) \we  (t_3t_4y_i) -
(t_1t_2y_i) \we  (t_3t_4x_i), \\ 
f_{H}(t)  & = &  
\sum_{i=1}^g 
(t_4t_1x_i) \we  (t_2t_3y_i) - 
(t_4t_1y_i) \we  (t_2t_3x_i), \\ 
f_X(t)  & = &  
\sum_{i=1}^g 
(t_1t_3x_i) \we  (t_4t_2y_i) - 
(t_1t_3y_i) \we  (t_4t_2x_i).
\end{eqnarray*}
It is easy to see that the map
$f_{IHX}:=f_I+f_H+f_X$ factors through $\Sym^4H\cong [4]_{\sp}$,
which is mapped onto the $[4]_{\sp}$-component of 
$\La^2\Sym^3H$.  
The proof of Theorem \ref{thm.new2}  goes exactly in the same way as
Theorem \ref{thm.new1}, where the only one point to be noted 
is that, for any embedding $I\hookrightarrow \Ga$, the naturally
induced total ordering $\tau_H$ on $\Ga_H$ is {\em not} 
$\Sym$-admissible, thus introducing the $IHX$ relation of figure \ref{ASIHX}.


\ifx\undefined\bysame
	\newcommand{\bysame}{\leavevmode\hbox 
to3em{\hrulefill}\,}
\fi


\begin{thebibliography}{[EMSS]}



\bibitem[AN]{AN} M.Asada, H.Nakamura,
     {\em On graded quotient modules of mapping class groups of surfaces},
     Israel J. Math. {\bf 90} (1995) 93--113.

\bibitem[B-N]{B-N} D. Bar-Natan, 
        {\em On the Vassiliev knot invariants}, Topology {\bf 34} 
        (1995)        423--472.

\bibitem[B]{B} B. Bollob\'{a}s,
        {\em Extremal Graph Theory}, Academic Press, 1978.

\bibitem[FH]{FH} W. Fulton, J. Harris,
        {\em Representation theory, a first course}, GTM {\bf 
        129},         Springer-Verlag, 1991.

\bibitem[GL1]{GL1} S. Garoufalidis and J. Levine,
        {\em Finite type 3-manifold invariants, the mapping 
        class group and blinks}, 
        Journal Diff. Geom. {\bf 47} (1997) 257--320.
        

\bibitem[GL2]{GL2} \bysame,
        {\em Finite type 3-manifold invariants and the Torelli group I},
        Inventiones, {\bf 131} (1998) 541--594.

\bibitem[GO]{GO} S. Garoufalidis, T. Ohtsuki,
       {\em On finite type 3-manifold invariants III: manifold weight 
       systems}, 
       Topology, {\bf 37} (1998) 227--244.

\bibitem[Ha]{Ha} R. Hain,
       {\em Infinitesimal presentations of the Torelli groups},
       Journal of AMS, {\bf 10} (1997) 597--651.

\bibitem[HL]{HL} R.Hain, E.Looijenga
	{\em Mapping class groups and moduli spaces of curves},
	Algebraic Geometry--Santa Cruz 1995. Proc. Sympos. Pure Math.
        {\bf 62} (1997) 97--142.

\bibitem[Ka]{Ka} M. Kapranov,
       {\em Rozansky-Witten invariants via Atiyah classes}, preprint
       March 1997.

\bibitem[K]{K} N.Kawazumi,
       {\em A generalization of the Morita-Mumford classes to
        extended maapping class groups for surfaces}, 
       Inventiones Math. {\bf 131} (1998) 137--149.

\bibitem[KM]{KM} N. Kawazumi, S. Morita,
       {\em The primary approximation to the cohomology of the 
       moduli space of curves and stable characteristic classes}, 
       Math. Research Letters, {\bf 3(5)} (1996), 629--642.

\bibitem[Ko1]{Ko1} M.~Kontsevich,
	{\em Vassiliev's knot invariants},
	Adv.\ in Sov.\ Math., {\bf 16(2)} (1993), 137--150.

\bibitem[Ko2]{Ko2} \bysame,
       {\em Formal (non)-commutative symplectic geometry},
       Gelfand Math. Seminars, 1990-92, Birkhauser, Boston, 
       (1993) 173-188.

\bibitem[Ko3]{Ko3} \bysame,
       {\em Feynmann diagrams and low-dimensional topology},
       Proceedings of the first European Congress of 
       Mathematicians,
       vol. 2, Progress in Math. {\bf 120} Birkhauser, 
       Boston, (1994)        97-121.

\bibitem[Ko4]{Ko4}\bysame,
       {\em Rozansky-Witten invariants via formal geometry},
       {\tt dg-ga/9704009}.

\bibitem[LMO]{LMO}  T.T.Q.  Le, J. Murakami, T. Ohtsuki,
       {\em A universal quantum invariant of 3-manifolds},
       Topology, {\bf 37} (1998) 539--574.

\bibitem[L]{L} T.T.Q.  Le, 
        {\em An invariant of \ihs s which is universal for all \fti s},
                Soliton Geometry and Topology: On the crossroad, 
        AMS Translations {\bf 2} Eds.
        V. Buchstaber, S. Novikov, 75--100.

\bibitem[Mo]{Mo} S. Morita,
       {\em A linear representation of the mapping class group of
       orientable surfaces and characteristic classes of vector bundles},
       Topology of Teichm\"{u}ller spaces,
       S. Kojima et al editors,  World Scientific,  (1996) 159-186.

\bibitem[Oh]{Oh} T. Ohtsuki,
        {\em Finite type invariants of integral homology 3-
        spheres}, 
        J. Knot Theory and its  Rami. {\bf 5} (1996) 101-115. 

\bibitem[RW]{RW} L. Rozansky, E. Witten,
        {\em Hyper-K\"{a}ehler geometry and invariants of 3-manifolds},
        Selecta Math. {\bf 3} (1997) 401--458.

\bibitem[W]{W} H.Weyl,
	{\em The Classical Groups; their invariants and representations},
	2nd ed., Princeton Univ. Press, Princeton NJ, 1946.

\bibitem[Wi]{Wi} E. Witten,
	{\em Quantum field theory and the Jones polynomial},
	Commun.\ Math.\ Phys.\ {\bf 121} (1989) 360--376.

 

\end{thebibliography}
\end{document}